\documentclass[12pt,draftcls,onecolumn]{IEEEtran}
\usepackage{cite}
\usepackage{times}
\usepackage{makeidx}
\usepackage{amsmath}
\usepackage{graphics}
\usepackage{graphicx}
\usepackage{amssymb}
\usepackage{euler}
\usepackage{supertabular}
\usepackage{epsfig}
\usepackage{delarray}
\usepackage{amstext}

\newcommand{\U}[1]{\,^{{#1}\hspace*{-1pt}}}

\title{\LARGE \bf Probabilistic Synchronization of Pinning Control}

\author{\authorblockN{Randa Herzallah}\\
\authorblockA{NCRG \\ Aston University, UK \\Email:r.herzallah@aston.ac.uk }
}

\begin{document}

\maketitle
\thispagestyle{empty}
\pagestyle{empty}

\section{Abstract}
This paper is concerned with synchronization of complex stochastic dynamical networks in the presence of noise and functional uncertainty. A probabilistic control method for adaptive synchronization is presented. All required probabilistic models of the network are assumed to be unknown therefore estimated to be dependent on the connectivity strength, the state and control values. Robustness of the probabilistic controller is proved via the Liapunov method. Furthermore, based on the residual error of the network states we introduce the definition of stochastic pinning controllability. A coupled map lattice with spatiotemporal chaos is taken as an example to illustrate all theoretical developments. The theoretical derivation is complemented by its validation on two representative examples.

\section{Introduction}\label{Sec:intro}

Current and future emergent systems in nature can be viewed as complex dynamical networks which are composed of nodes representing individuals in the system and links representing interactions among them. As a typical collective behavior in nature, synchronization has been one of the key issues that is investigated in the literature. Over the past decade chaos control, stability (controllability)~\cite{Jun11,Jinhu05}, and synchronization~\cite{Wenwu09,Xuan11,Yu08,Zhi07} of large scale dynamic systems have been some of the issues with wide interest due to their potential applications in power grid systems, biological systems, reaction diffusion systems and many others. Typically because of the large number of nodes in real world complex dynamical networks, it has been widely believed that it is impossible to add controllers to all nodes. Consequently, pinning control has been proposed as a viable strategy which requires a significantly smaller number of controllers which are injected to a fraction of network nodes. For a more detailed discussion on pinning control see~\cite{Sorrentino07,Xiang07,Tang09,Murizio08,Abhijit10,Wenlian10,Chen07}.

In particular, previous pinning control methodologies are based mainly on linear feedback control theory which is served as a simple and effective approach for stabilization and synchronization. In~\cite{Xiang07}, the stabilization problem of complex dynamical networks with general coupling topology by pinning a small fraction of nodes with local negative feedback controllers was discussed. The result that a network under a typical framework can realize synchronization subject to any linear feedback pinning scheme by using adaptive tuning of the coupling strength was proved in~\cite{Yu09}. These studies largely focus on developing deterministic control algorithms, therefore, necessitate perfect information exchange among the nodes. Since information exchange between the nodes in the network are normally affected by noise and quantization error, the development of control algorithms that can handle uncertain information has recently received great deal of attention in the field. Examples include the extended result to the density of pinning sites for spatially chaotic systems with linear quadratic optimal control, provided in~\cite{Grigoriev97}. It is shown that the minimal density of the pinning nodes depends on the strength of noise in the system. An adaptive control approach by adding an additional control input to all or a fraction of nodes in the network was proposed in~\cite{Zhou06} to guarantee synchronization in the presence of uncertainties.  Another control scheme for synchronization of complex networks was proposed in~\cite{DeLellis09}. It is based on the local adaptive control approach. This scheme is claimed to be useful in applications where the coupling gains cannot be chosen globally, or where synchronization needs to be attained in the presence of uncertainties and noise~\cite{DeLellis09}.

However, in all of the control strategies described above either the adaptive control method was implemented to solve the synchronization problem of uncertain systems or only uniform noise was assumed to be acting on the dynamics of the systems. Moreover, the performance objectives concerned were confined to be the mean value of the stochastic output and consequently control is still confined to the deterministic control strategy. In this paper, by using neural network methods and a cost function derived from a Kullback--Leibler distance between the joint probability density function of the closed-loop system and an ideal joint probability density function, a probabilistic controller is proposed for a stochastic uncertain complex dynamical network. This method was originally proposed in~\cite{Herzallah11} to obtain a general solution for stochastic systems subject to random inputs and deterministic systems characterized by functional uncertainty with unknown probability density functions. The synchronization problem of complex dynamical networks, on the other hand, was not discussed.

The consensus problem of complex networks in the presence of formation leaders is shown in~\cite{Jadbabaie03,Lafferriere05} to be similar to the pinning control problem. The formation leader defines the common value that all other nodes in the network should converge to. Consensus problems in a stochastic settings have received more attention. This can be referred to that the network dynamics in consensus problems is linear while in pinning control the coupled systems can be strongly nonlinear which complicates the problem further. In~\cite{Huang2010} a stochastic approximation type algorithms are employed for individual states to converge in mean square to the same limit. The agreement problem over random information networks, where the existence of information channel between a pair of units is probabilistic is discussed in~\cite{Hatano05} and the asymptotic agreement of the network is addressed via notions of stochastic stability. Necessary and sufficient conditions for mean square consentability of the averaging protocol for a stochastic directed network is discussed in~\cite{Abaid2011}. The distributed average consensus for random topologies and noisy channels is studied in~\cite{Kar09}, where noisy consensus is shown to lead to a bias--variance dilemma. It is worth noticing that these researches are mainly concerned with achieving consensus by designing the optimal link weights when the connectivity graph of the network is fixed or random. The current work however, is concerned with achieving consensus by designing a probabilistic pinning control algorithm.

Several research efforts have been devoted to the analysis of synchronization of complex networks under uncertain and noisy conditions by adding additional deterministic control input to all or fraction of nodes in the network~\cite{Grigoriev97,Zhou06,DeLellis09}, as discussed earlier. More often, under noisy and uncertain conditions, it is more appropriate and effective to consider probabilistic rather that deterministic control algorithms. To the best of the authors knowledge, this has not been considered in earlier work in the field of pinning controlling complex dynamical networks.

To summarize, the aim of this paper is to study the problem of how to achieve synchronization of complex stochastic dynamical networks via design of probabilistic controller, expanding and integrating the results presented in~\cite{Herzallah11}. In particular, the stochastic controllability condition is introduced and the stability of the proposed probabilistic control methodology is proven. Compared with the existing results on the topic, this paper has three distinct features that have not been reported in literature. Firstly, based on the theory of dynamic programming and the Liapunov method  we design fully probabilistic controllers, as opposed to current existing deterministic controllers, that are applied to a subset of nodes in the network. It will be demonstrated that these probabilistic controllers guarantee synchronization in the presence of noise and functional uncertainties. Secondly, the probabilistic models of the forward stochastic dynamics of the network are assumed to be unknown and hence estimated adaptively as a function of the connectivity strength, state values, and pinning control nodes. Thirdly, the Kullback–-Leibler distance between the joint probability density function of the closed-loop system and an ideal joint probability density function is used as the cost function to relax the conservativeness on the assumed generative distribution of the random noise. The proposed probabilistic control methodology will be applied to controlling spatiotemporally chaotic systems with the coupled map lattice taken as an example. To emphasize, probabilistic controllers proposed in this paper provide a pragmatic method for achieving synchronization in complex stochastic dynamical networks under uncertain conditions. It will be seen, that despite of the uncertainty embedded in the information states of those stochastic networks at each instant of time one can still prove synchronization results that resemble those for the certain and deterministic case.

The rest of the paper is organized as follows. In Section~\ref{Sec:ProbForm} the model description and some preliminaries are give. The pinning controllability of stochastic networks is introduced in Section~\ref{Sec:PinCont}. Synchronization of probabilistic control problem is given in Section~\ref{Sec:SynchProbCont}. In Section~\ref{Sec:AsympStab} the proof of asymptotic stability of the proposed probabilistic controller is detailed. Section~\ref{Sec:Results} contains some
simulation results to show the effectiveness of the proposed controller. The conclusion is provided in Section~\ref{sec:Con}.

\section{Model Description and Preliminaries}\label{Sec:ProbForm}
Consider the coupled map lattice consisting of $L$ lattice nodes with periodic boundary conditions~\cite{Grigoriev97,Gang94,Kaneko84},
\begin{eqnarray}
z^i_{t+1} &=&  F(z^{i-1}_t,z^i_t,z^{i+1}_t) \nonumber \\ &=& f[(1-2 \epsilon) z^i_t + \epsilon(z^{i-1}_t + z^{i+1}_t)],
\label{eq:eq1}
\end{eqnarray}
where $i=1,2,\dots,L$ are the lattice sites, $L$ is the system size, and the periodic boundary conditions are given by $z^{i+L}_t=z^i_t$. The local map $f(z)$ is defined to be a nonlinear function of the following form
\begin{equation}
f(z) = a z (1-z).
\label{eq:eq2}
\end{equation}
This coupled map lattice exhibits chaotic characteristics in the regime $3.57 < a \le 4.0$ and has a homogeneous steady state $z^\star = 1 - 1/a$. The aim is to synchronize network~\eqref{eq:eq1} onto a homogeneous stationary state such that \begin{equation}
\lim_{t \rightarrow \infty} \parallel z^i_t - z^\star \parallel = 0,~~~~~~~~~~~~~~~~~i=1,2,\dots, L.
\end{equation}
The first attempt in this direction was reported in~\cite{Gang94}, where $M$ periodically pinning control actions are applied at sites $\{i_1,\dots,i_M\}$ in the following way
\begin{eqnarray}
z^i_{t+1} = F(z^{i-1}_t,z^i_t,z^{i+1}_t) + \sum_{m=1}^M \delta(i-i_m) u^m_t,
\label{eq:eq3}
\end{eqnarray}
where $u^m$ is the control action applied at site $m$. However, this required a very dense array of pinning controllers $L_p = L/M$ in the physically interesting interval of parameters $3.57 < a \le 4.0$. To further understand how pinning should be placed, a linearized form of~\eqref{eq:eq1}, about the homogeneous steady state ${\bf z}_t= (z^{\star^1},\dots,z^{\star^L})$, is considered in~\cite{Grigoriev97} and the following standard linear equation is obtained:
\begin{equation}
{\bf x}_{t+1} = \tilde{A} {\bf x}_t + \tilde{B} {\bf u}_t,
\label{eq:eq4}
\end{equation}
in which ${\bf x}={\bf z}-{\bf z}^\star$ represents the state vector, the $L\times L$ Jacobian matrix $\tilde{A}$, is given by
\begin{eqnarray}
\tilde{A} = \alpha \left [ \begin{array}{ccccc} 1-2\epsilon & \epsilon & 0 &\dots & \epsilon \\ \epsilon & 1-2 \epsilon & \epsilon &\dots & 0 \\ 0 & \epsilon & 1-2\epsilon & \dots & 0 \\ \vdots & \vdots & \vdots & \ddots & \vdots \\ \epsilon & 0 & 0 & \dots &1-2\epsilon
\end{array} \right ],
\label{eq:eq5}
\end{eqnarray}
where $\alpha= \frac{\partial f(z)}{\partial z} \mid_{z=z^\star}$, and $\tilde{B}$ is an $L\times M$ control matrix with $\tilde{B}_{ij} = \sum_m \delta(j-m) \delta (i-i_m)$ which depends on how the pinning sites are placed. Thus, to minimize the number of pinning nodes, the limits of the control scheme is extended in~\cite{Grigoriev97} by making the system controllable as opposed to stabilizable. The controllability condition of the system~\eqref{eq:eq4} is that the rank of the following controllability matrix equal to L, or
\begin{equation}
\text{rank}(C[\tilde{A},\tilde{B}]) = \text{rank}[\tilde{B} \vdots \tilde{A}\tilde{B} \vdots \cdots \vdots \tilde{A}^{L-1}\tilde{B}] = L.
\end{equation}
Since the system~\eqref{eq:eq1} has parity symmetry, it is shown~\cite{Grigoriev97} that the controllability condition is satisfied using at least two pinning sites, $M=2$, with the only restriction that $L$ should not be multiple of $|i_2-i_1|$; or otherwise the mode with the period $2|i_2-i_1|$ becomes uncontrollable.

On the other hand, the existence of an external random input that affects the dynamics of the systems is very common in many areas of science and engineering. Therefore our aim in this paper is to further discuss the controllability condition of the stochastic version of~\eqref{eq:eq4} and propose a more general fully probabilistic control algorithm.

\section{Pinning Controllability of Stochastic Networks}\label{Sec:PinCont}
Further to the approach presented in~\cite{Grigoriev97} we define and assess the stochastic pinning controllability of the coupled map lattice networks defined in~\eqref{eq:eq1} with additive external random inputs. For that purpose we consider the general time varying dynamics of the linearized equation~\eqref{eq:eq4} with external random inputs added to its right hand side as follows:
\begin{equation}
{\bf x}_{t+1} = \tilde{A}_t {\bf x}_t + \tilde{B}_t {\bf u}_t+ \tilde{E}_t \tilde{\kappa}_{t},
\label{eq:eq4Stoch}
\end{equation}
where $\tilde{A}, \tilde{B}, {\bf x}$ and ${\bf u}$ have similar definitions as before, $\tilde{\kappa}_{t}$ is an additive noise signal assumed to have zero mean Gaussian distribution of covariance $\tilde{\Sigma}$, and $\tilde{E}$ is the noise matrix . The solution to~\eqref{eq:eq4Stoch} can be easily verified to be given by
\begin{eqnarray}
{\bf x}_{t+L}&= \phi(t+L,t) {\bf x}_t+\sum_{j=0}^{L-1} \phi(t+L,j+t+1) \tilde{E}_{j+t}\tilde{\kappa}_{j+t}
\nonumber \\ &  + \sum_{j=0}^{L-1} \phi(t+L,j+t+1) \tilde{B}_{j+t}{\bf u}_{j+t},
\label{eq:StochSol}
\end{eqnarray}
with
\begin{equation}
\phi(t+L,t)=\tilde{A}_{t+L-1}\tilde{A}_{t+L-2}....\tilde{A}_{t}.
\end{equation}
For the time invariant system where $\tilde{A}$, and $\tilde{B}$ are constant matrices we have
\begin{equation}
\phi(t+L,t)=\tilde{A}^L.
\end{equation}
Because of the noise disturbances the state vector at time $t+L$ differs from $\phi(t+L,t){\bf x}_t$, the state vector at time $t+L$ with no disturbances by
\begin{eqnarray}
d_L &=& {\bf x}_{t+L} - \phi(t+L,t){\bf x}_t, \nonumber \\
   &=&\sum_{j=0}^{L-1} \phi(t+L,j+t+1) \tilde{E}_{j+t}\tilde{\kappa}_{j+t}.
\label{eq:ResErr}
\end{eqnarray}
Because $\tilde{\kappa}_{t}$ is an additive Gaussian noise signal with zero mean, the random variable $d_L$ is also Gaussian with
\begin{eqnarray}
E(d_L) &=& 0, \nonumber \\
\text{Cov}(d_L) &=& \sum_{j=t+1}^{t+L} \phi(t+L,j) \tilde{E}_{j-1} \tilde{\Sigma}_{j-1} \tilde{E}_{j-1}^{T}\phi^T(t+L,j). \nonumber
\end{eqnarray}
Note that, obviously, for complete controllability of the stochastic network~\eqref{eq:eq4Stoch}, the residual error defined in~\eqref{eq:ResErr} should remain bounded. Within this framework, we can introduce the following mathematical definition for the concept of stochastic controllability.

\textit{\textbf{Definition $1$:}} Let $\text{Cov}(d_L)=\Psi_L$, then if the matrix $\Psi_L$ is positive definite and $\parallel \Psi_L \parallel$ remains bounded for all $L$, where $\parallel.\parallel$ is a norm in an Euclidean space, then the residual error $d_L$ will remain bounded for all $L$~\cite{Masanao67}. Therefore $\Psi_L$ is called the stochastic controllability matrix of the stochastic dynamical network~\eqref{eq:eq4Stoch}.

For the stationary system where $\tilde{A}, \tilde{B}, \phi$ and $\tilde{E}$ are constant matrices and $\tilde{\Sigma}>0$ the criterion of complete controllability is that the rank of the stochastic controllability matrix, $SC[\phi,\tilde{E}]$ equal to $L$,
\begin{equation}
\text{rank}(SC[\phi,\tilde{E}])= L.
\end{equation}

\section{Synchronization Probabilistic Control Problem}\label{Sec:SynchProbCont}
Consider a stationary dynamics of the time variant stochastic coupled map lattice defined in~\eqref{eq:eq4Stoch} as,
\begin{equation}
{\bf x}_{t+1} = \tilde{A} {\bf x}_t + \tilde{B} {\bf u}_t+ \tilde{\kappa}_{t},
\label{eq:eq4Stoch1}
\end{equation}
where here without loss of generality $\tilde{E}$ is taken to be the identity matrix. The synchronization control problem confronted here is to design a control strategy for pinning nodes to synchronize the state of the coupled map lattice to the homogenous state, i.e $z\rightarrow z^*$ or $ {\bf x} \rightarrow 0$.  Then according to Definition $1$ all nodes should synchronize such that the residual error, $d_L$ of the stochastic dynamics is bounded. However, because of the noise input $\tilde{\kappa}_{t}$ the present state and present and future controls do not completely specify the future state, but instead determine only the probability distribution of these states, $s({\bf x}_{t+1} \mid {\bf u}_t, {\bf x}_t)$ . It is assumed that the node disturbances $\tilde{\kappa}_{t}$ are unknown and hence the probability distribution of the states is unknown. Thus, the synchronization strategy must be robust to unmodelled stochastic dynamics and unknown disturbances.

The gaol of this synchronization control problem fits naturally with the probabilistic control theory. Therefore, to achieve this goal we consider designing a probabilistic controller $c({\bf u}_t \mid {\bf x}_t)$ such that the joint probability density function (pdf) of the closed loop system, $f(D) = \prod_{t=0}^{L} s({\bf x}_{t+1} \mid {\bf u}_t, {\bf x}_t) c({\bf u}_t \mid {\bf x}_t)$ is made as close as possible to a desired pdf, $\U{I}f(D) = \prod_{t=0}^{L} \U{I}s({\bf x}_{t+1} \mid {\bf u}_t, {\bf x}_t)\U{I}c({\bf u}_t \mid {\bf x}_t)$. This design method was originally presented in~\cite{Karny96,Herzallah11}, where the probabilistic controller is obtained such that it minimizes the following cost to go function derived from the Kullback-Leibler divergence (KLD) between the actual joint pdf $f(D)$ and the ideal joint pdf $\U{I}f(D)$ and defined as follows,
\begin{eqnarray*}
&&-\ln(\gamma({\bf x}_{t}))=\min_{\left\{c({\bf u}_t \mid {\bf x}_t)_{t\geq
\tau}^{L} \right\}} \sum_{t=\tau}^{L}\int
f({\bf x}_{t+1}, {\bf u}_t, \dots, {\bf x}_{L+1}, {\bf u}_L) \\
&&\times\ln\left(
\frac{s({\bf x}_{t+1} \mid {\bf u}_t, {\bf x}_t) c({\bf u}_t \mid {\bf x}_t)}
     {\U{I}s({\bf x}_{t+1} \mid {\bf u}_t, {\bf x}_t)\U{I}c({\bf u}_t \mid {\bf x}_t)}
\right)\,\mathrm{d}({\bf x}_{t+1}, {\bf u}_t, \dots, {\bf x}_{L+1}, {\bf u}_L), \\ &&
\end{eqnarray*}
for arbitrary $t \in \{1, \dots, L\}$. Using this definition minimization is then performed recursively to give the following recurrence functional equation~\cite{Herzallah11}
\begin{eqnarray}
\label{eq:probcritic4} &&-\ln(\gamma({\bf x}_{t}))=
\min_{c({\bf u}_t \mid {\bf x}_t)} \int
s({\bf x}_{t+1} \mid {\bf u}_t, {\bf x}_t) c({\bf u}_t \mid {\bf x}_t) \nonumber \\
&&\times\bigg[\underbrace{\ln\left(
\frac{s({\bf x}_{t+1} \mid {\bf u}_t, {\bf x}_t) c({\bf u}_t \mid {\bf x}_t)}
     {\U{I}s({\bf x}_{t+1} \mid {\bf u}_t, {\bf x}_t)\U{I}c({\bf u}_t \mid {\bf x}_t)}
\right)}_{\equiv \mbox{partial cost $\Longrightarrow U({\bf u}_t, {\bf x}_t)$}}
-\underbrace{\ln(\gamma({\bf x}_{t+1}))}_{\mbox{optimal
cost-to-go}}\bigg] \,\mathrm{d}({\bf x}_{t+1},{\bf u}_{t}).
\end{eqnarray}
Furthermore, to overcome the curse of dimensionality problem arising from the probabilistic control~\eqref{eq:probcritic4}, adaptive critic methods were proposed in~\cite{Herzallah11} to approximate the optimal cost to go function and the probabilistic controller. Unknown pdfs were also estimated using recent development from neural network models. Numerical experiments and previous analytical studies have shown the usefulness of this control approach to obtain the control efforts for general stochastic systems subject to functional uncertainty and random inputs~\cite{Herzallah11}.  The synchronization problem of complex dynamical systems on the other hand was not discussed. Furthermore, A pressing open problem to the probabilistic control approach is to prove analytically the asymptotic stability of the controlled system.

We start by first discussing the estimation problem of unknown probabilistic models of the lattice network and giving some definitions that will be used throughout the rest of the paper. To estimate the probabilistic model of the lattice network we adopt the method proposed in~\cite{Herzallah11}, where neural network models are used to provide a prediction for the conditional expectation of the system output and calculating the global average variance of its residual error. For the linearized coupled map lattice~\eqref{eq:eq4Stoch1} generalized linear neural network (GLNN) models can be used to provide estimates for all unknown pdfs. Therefore using GLNN models, the stochastic model of the coupled map lattice defined in~\eqref{eq:eq4Stoch1} can be shown to be given by,
\begin{eqnarray}
{\bf x}_{t+1} &=& A {\bf x}_{t} + B {\bf u}_{t} +  \kappa_{t}, \label{eq:probcritic7}
\end{eqnarray}
where $A {\bf x}_{t} + B {\bf u}_{t}$ represents estimated conditional expectation of the lattice state with $A$ and $B$ being estimate of the jacobian and control matrices respectively, $\kappa_{t}$ is the noise of the residual error of the lattice output which is shown~\cite{Herzallah11} to be close to Gaussian random noise
with zero mean and $\Sigma$ covariance. According to the developed stochastic model in~\eqref{eq:probcritic7}, the distribution of the state values will be Gaussian distribution with expected means provided by the approximating GLNN and a global covariance $\Sigma$ given by the residual value of the error between actual and estimated
states $E[\{{\bf x}_{t+1} -(A {\bf x}_{t} + B {\bf u}_{t})\}\{{\bf x}_{t+1} -(A {\bf x}_{t} + B {\bf u}_{t})\}^T]$,
\begin{equation}
s({\bf x}_{t+1} \mid {\bf x}_{t},{\bf u}_{t})  \rightsquigarrow \mathcal{N}_{{\bf x}_{t+1}}
(A{\bf x}_t+ B {\bf u}_{t}, \Sigma). \label{eq:probcritic7p}
\end{equation}

For the synchronization problem of the coupled map lattice, the ideal state distribution is assumed to be given by
\begin{equation}
\U{I}s({\bf x}_{t+1} \mid {\bf u}_t, {\bf x}_t) \rightsquigarrow \mathcal{N}_{{\bf x}_{t+1}}(0, \Sigma). \label{eq:idealForDyn}
\end{equation}
It reflects the synchronization problem with the aim of reaching the zero state, with the spread determined by the covariance of innovations $\Sigma$. The randomized controller to be designed is described by the following stochastic model
\begin{eqnarray}
{\bf u}_{t} &=& C {\bf x}_{t} + \omega_{t}, \nonumber \\ c({\bf u}_{t}|{\bf x}_{t}) & \rightsquigarrow & \mathcal{N}_{{\bf u}_{t}}(C {\bf x}_{t}, \Gamma),
\label{eq:probcritic11}
\end{eqnarray}
where C is the matrix of the controller parameters, $\omega_{t}$ is the residual error of the control input, and $\Gamma$ is the covariance of the residual error of control. The distribution of the ideal controller is also
assumed to be
\begin{equation}
\U{I}c({\bf u}_{t}|{\bf x}_{t}) \rightsquigarrow \mathcal{N}_{{\bf u}_{t}} (0, \Gamma). \label{eq:idealConDist}
\end{equation}

\section{Proof of Asymptotic Stability}\label{Sec:AsympStab}
 If we use the randomized controller~\eqref{eq:probcritic11}, the governing stochastic equation of the network~\eqref{eq:probcritic7} can be recast as follows:
\begin{equation}
{\bf x}_{t+1} = (A + B C) {\bf x}_{t} + B \omega_{t} + \kappa_{t}.
\label{eq:eqproof3}
\end{equation}
The probabilistic control design problem is then to determine a stabilizable probabilistic model $c({\bf u}_{t}|{\bf x}_{t})$ of the randomized controller under the following assumption:

\textit{\textbf{Assumption $1$:}} All eigenvalues, $\lambda_k$ of the matrix $(A + BC)$ in~\eqref{eq:eqproof3} lie inside the unit circle. Similar to the asymptotic stability of deterministic control systems~\cite{Ogata87} this is equivalent to,
\begin{equation}
\text{rank}(AS[((\Sigma^{-1})^{1/2})^T, A]) = L.
\end{equation}

The above control objective can be achieved by minimization of~\eqref{eq:probcritic4} subject to the constraint equations specified by~\eqref{eq:probcritic7p}, \eqref{eq:idealForDyn} and~\eqref{eq:idealConDist}. This, leads to the probabilistic optimal feedback control law specified in the following theorem.

\textit{\textbf{Theorem $1$:}} The feedback control law minimizing optimal cost to go function~\eqref{eq:probcritic4} subject to the pdf of the system dynamics~\eqref{eq:probcritic7p} and ideal distributions given by~\eqref{eq:idealForDyn} and~\eqref{eq:idealConDist}
\begin{eqnarray}
{\bf u}_{t} &=& C {\bf x}_{t} + \omega_{t},
\label{eq:probcritic11T}
\end{eqnarray}
with
\begin{eqnarray}
C &=& -(B^T M B + B^T \Sigma^{-1} B + \Gamma^{-1})^{-1} (B^T M A + B^T \Sigma^{-1} A), \label{eq:ContPar1}
\end{eqnarray}
is a stabilizing control law and
\begin{eqnarray}
-\ln(\gamma({\bf x}_{t})) &=& 0.5 {\bf x}_{t}^T M {\bf x}_{t} + \mathbb{Q}_0,
\label{eq:probcritic15}
\end{eqnarray}
with
\begin{equation}
M = A^T \Sigma^{-1} A + A^T M A - (A^T M B + A^T \Sigma^{-1} B) (B^T M B + B^T \Sigma^{-1} B + \Gamma^{-1})^{-1} (B^T M A + B^T \Sigma^{-1} A), \label{eq:QuadMat}
\end{equation}
is the quadratic cost function. Here $\mathbb{Q}_0 \ge 0$ is some positive constant.

\textit{\textbf{Proof:}} To prove the theorem we start by evaluating the partial cost $U({\bf u}_t, {\bf x}_t)$ in~\eqref{eq:probcritic4} repeated here,
\begin{eqnarray}
U({\bf u}_t, {\bf x}_t) = \int s({\bf x}_{t+1} \mid {\bf u}_t, {\bf x}_t) c({\bf u}_t \mid {\bf x}_t)\times \ln\left(
\frac{s({\bf x}_{t+1} \mid {\bf u}_t, {\bf x}_t) c({\bf u}_t \mid {\bf x}_t)}
     {\U{I}s({\bf x}_{t+1} \mid {\bf u}_t, {\bf x}_t)\U{I}c({\bf u}_t \mid {\bf x}_t)}
\right) \mathrm{d}({\bf x}_{t+1},{\bf u}_{t}). \label{eq:eqproof2}
\end{eqnarray}
The Evaluation of~\eqref{eq:eqproof2} yields,
\begin{equation}
U({\bf u}_t, {\bf x}_t) = {\bf x}_t^T [C^T \Gamma^{-1} C + (A + BC)^T \Sigma^{-1} (A+BC)] {\bf x}_t, \label{eq:eqproof4}
\end{equation}
where we made use of~\eqref{eq:probcritic11T} and where $\Gamma$ and $\Sigma$ are positive definite real symmetric matrices. Under Assumption $1$, there exists a Liapunov function that is positive definite and whose derivative is negative definite as follows,
\begin{eqnarray}
{\bf x}_t^T [C^T \Gamma^{-1} C + (A + BC)^T \Sigma^{-1} (A+BC)] {\bf x}_t &=  - \int[{\bf x}^T_{t+1} M {\bf x}_{t+1} - {\bf x}^T_t M {\bf x}_t] \mathrm{d}({\bf x}_{t+1},{\bf u}_{t}) \nonumber \\ & = {\bf x}_t^T[-(A + BC)^T M (A+BC) + M] {\bf x}_t. \label{eq:eqproof5}
\end{eqnarray}
Comparing the two sides of~\eqref{eq:eqproof5} and noting that this equation must hold for any ${\bf x}_t$, it is then required that
\begin{equation}
C^T \Gamma^{-1} C + (A + BC)^T \Sigma^{-1} (A+BC) = -(A + BC)^T M (A+BC) + M. \label{eq:eqproof6}
\end{equation}
Equation~\eqref{eq:eqproof6} can be modified as follows
\begin{eqnarray}
&A^T \Sigma^{-1} A + A^T M A - M + C^T (\Gamma^{-1} + B^T M B + B^T \Sigma^{-1} B) C  + C^T B^T \Sigma^{-1} A \nonumber \\ &+ A^T \Sigma^{-1} B C + C^T B^T M A + A^T M B C = 0. \label{eq:eqproof7}
\end{eqnarray}
This last equation can further be modified as follows
\begin{eqnarray}
&\bigg\{(\Gamma^{-1} + B^T M B + B^T \Sigma^{-1} B)^{1/2} C + (\Gamma^{-1} + B^T M B + B^T \Sigma^{-1} B)^{-1/2} (B^T M A + B^T \Sigma^{-1} A) \bigg\}^T \nonumber \\ &\times \bigg\{(\Gamma^{-1} + B^T M B + B^T \Sigma^{-1} B)^{1/2} C + (\Gamma^{-1} + B^T M B + B^T \Sigma^{-1} B)^{-1/2} (B^T M A + B^T \Sigma^{-1} A) \bigg\} \nonumber \\ &- (B^T M A + B^T \Sigma^{-1} A)^T (\Gamma^{-1} + B^T M B + B^T \Sigma^{-1} B)^{-1}(B^T M A + B^T \Sigma^{-1} A) \nonumber \\ &+ A^T \Sigma^{-1} A + A^T M A - M    = 0.  \label{eq:eqproof8}
\end{eqnarray}
Minimization of the performance index~\eqref{eq:probcritic4} with respect to $C$ requires minimization of the left hand side of~\eqref{eq:eqproof8} with respect to $C$. Since
\begin{eqnarray}
&\bigg\{(\Gamma^{-1} + B^T M B + B^T \Sigma^{-1} B)^{1/2} C + (\Gamma^{-1} + B^T M B + B^T \Sigma^{-1} B)^{-1/2} (B^T M A + B^T \Sigma^{-1} A) \bigg\}^T \nonumber \\ &\times \bigg\{(\Gamma^{-1} + B^T M B + B^T \Sigma^{-1} B)^{1/2} C + (\Gamma^{-1} + B^T M B + B^T \Sigma^{-1} B)^{-1/2} (B^T M A + B^T \Sigma^{-1} A) \bigg\}, \label{eq:eqproof9}
\end{eqnarray}
is nonnegative, the minimum occurs when it is zero or when
\begin{eqnarray}
(\Gamma^{-1} + B^T M B + B^T \Sigma^{-1} B)^{1/2} C = -(\Gamma^{-1} + B^T M B + B^T \Sigma^{-1} B)^{-1/2} (B^T M A + B^T \Sigma^{-1} A).  \label{eq:eqproof10}
\end{eqnarray}
Solving this last equation for $C$ yields the optimal control law specified in~\eqref{eq:ContPar1}. Substitution of~\eqref{eq:ContPar1} into equation~\eqref{eq:eqproof7} gives~\eqref{eq:QuadMat}.

The proof is thus completed.

\section{Numerical Results}\label{Sec:Results}
This section will demonstrate the effectiveness of the probabilistic control approach specified by Theorem $1$ in the non chaotic and chaotic regimes of the lattice. It is shown that the approach effectively enforces synchronization using only two pinning nodes for unknown system dynamics and unknown random input at each node.

\subsection{Non Chaotic Map Lattice}
\begin{figure}[htbp]
  \centerline{
    \begin{tabular}{cc}
     \resizebox{3in}{2in}{\includegraphics{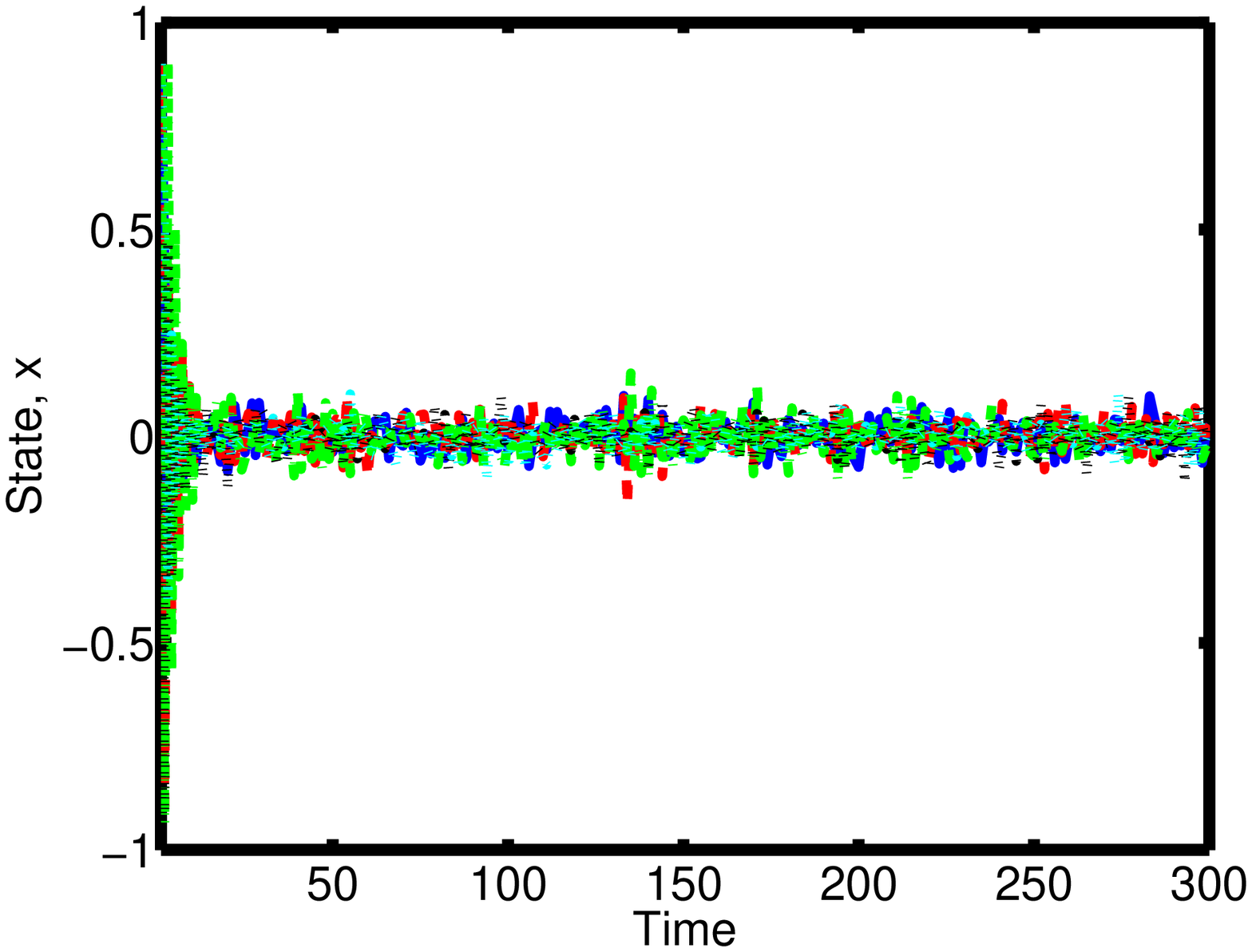}}
     &
      \resizebox{3in}{2in}{\includegraphics{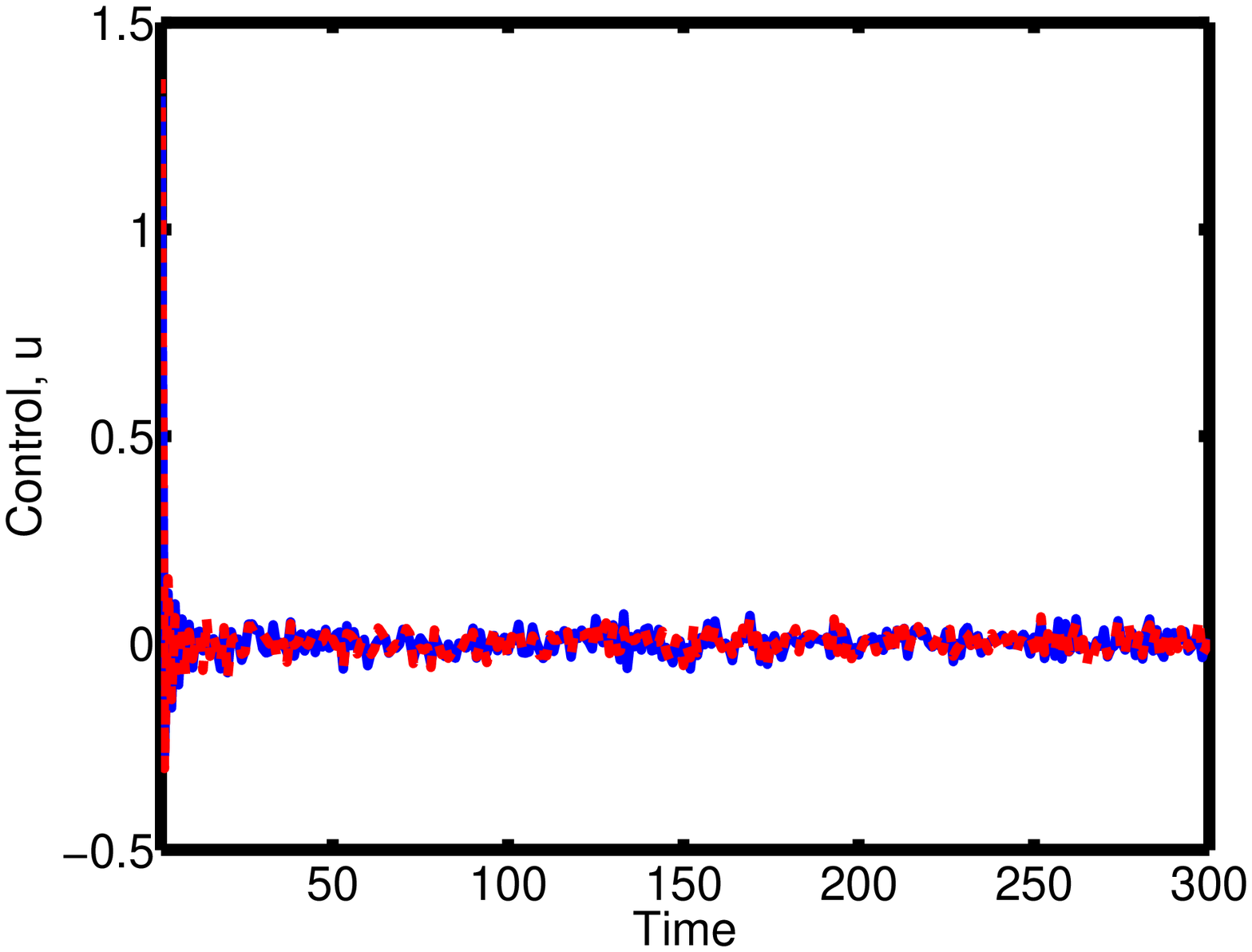}}
      \\
      (a) & (b)
      \\
     \resizebox{3in}{2in}{\includegraphics{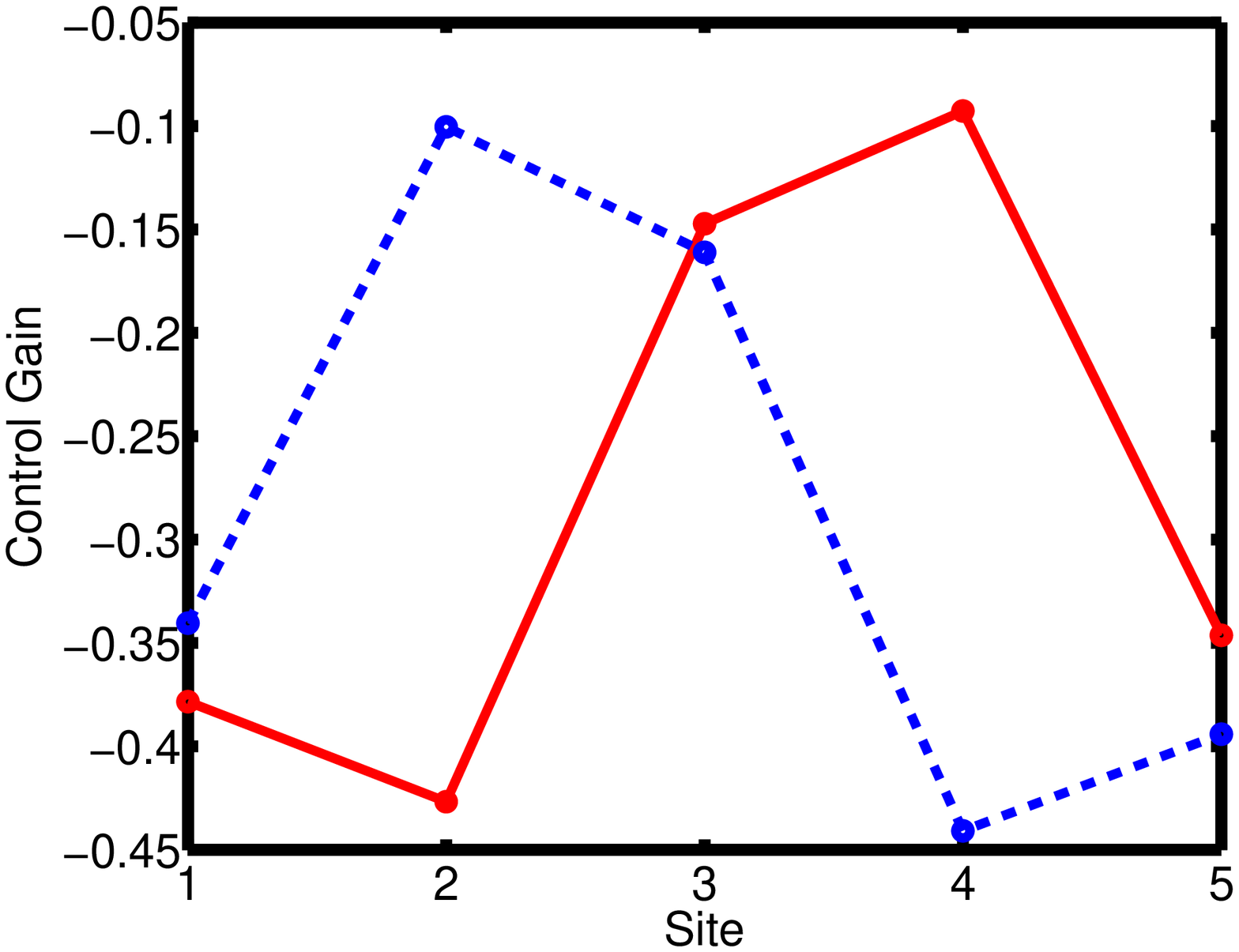}}
     &
     \resizebox{3in}{2in}{\includegraphics{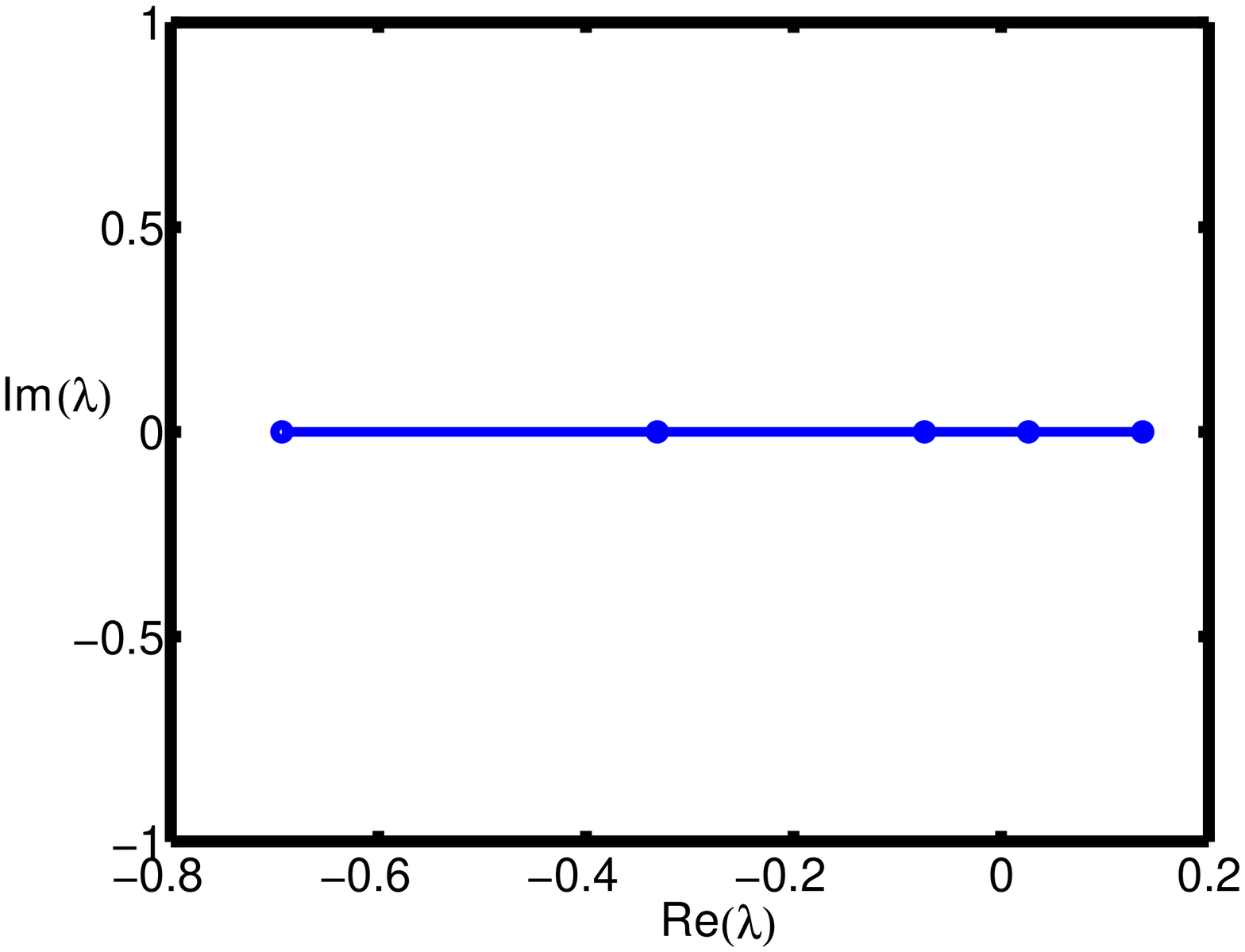}}
     \\
      (c) & (d)
     \end{tabular}
    }
  \caption[Non chaotic coupled map lattice.] {Non chaotic coupled map lattice with, $L=5$, $a=3$, and $\epsilon=0.33$: (a) Evolution of states. (b) Control efforts. (c) Control Gain. (d) Eigenvalues of the asymptotic stability matrix.}
  \label{fig:Fig1}
\end{figure}

The example considered here is for the logistic coupled map lattice, $f(z) = a z(1-z)$ in its non--chaotic regime with $a=3.0$, $\epsilon=0.33$ and $L=5$ and with an external Gaussian random input affecting the dynamics. The two control actions are placed next to each other at the sides of the lattice. Hence the equation of the coupled map lattice becomes:
\begin{equation}
{\bf x}_{t+1} = A {\bf x}_t + B {\bf u}_t + \kappa_{t},
\label{eq:eq19}
\end{equation}
where
\begin{eqnarray}
A = \left [ \begin{array}{ccccc} -0.34 & -0.33 & 0 &0 & -0.33 \\ -0.33 & -0.34 & -0.33 &0 & 0 \\ 0 & -0.33 & -0.34 & -0.33 & 0 \\ 0 & 0 & -0.33 & -0.34 & -0.33 \\ -0.33 & 0 & 0 & -0.33 &-0.34
\end{array} \right ], \nonumber \\ B = \left [ \begin{array}{cc} 1 &~0 \\ 0 &~0 \\ 0 &~ 0 \\ 0 &~ 0 \\ 0 &~1 \end{array} \right ]~~~~E[\kappa_{t} \kappa_{t}^T] = 0.001 I_{5 \times 5}.~~~~~~~~ \nonumber
\end{eqnarray}
The lattice is initially at time $t=0$ in state ${\bf x} = 0.9$ and the aim is to return the lattice state to the origin (the fixed point position) or a state close to the origin. As a first step in the solution, the Gaussian probability density function of the stochastic lattice model described by~\eqref{eq:eq19} is estimated off--line using generalized linear model and a global diagonal covariance matrix $\Sigma = [0.00095, 0.00105, 0.00097, 0.0009, 0.0011]$  as discussed in Section~\ref{Sec:SynchProbCont}. The covariance matrix of the controller is taken to be $\Gamma = 0.01 I_{2 \times 2}$. The states of the lattice network and the control efforts are illustrated in Figures~\ref{fig:Fig1}(a) and (b), respectively, which show that the controlled network is globally synchronized and that the designed probabilistic control approach is very effective.

Using the estimated global average covariance matrix, $\Sigma$ the rank of the asymptotic stability matrix, $\text{rank}(AS[((\Sigma^{-1})^{1/2})^T, A])$ is found to be $5$, hence satisfying Assumption$1$. This is further verified by finding the eigenvalues, $\lambda_k$ of the matrix $(A +BC)$ which all are found to be inside the unit circle as can be seen from Figure~\ref{fig:Fig1}(d). Figure~\ref{fig:Fig1}(c) gives the solution of the control matrix $C$. The control effort $C^{mi}$ is larger for those sites that are far away from the pinning site $i_m$. This is expected since the probabilistic feedback effort is applied indirectly through coupling to the neighbors which in turn means that the perturbation introduced by the controllers decays with the increasing distance to the pinning sites. This agrees with the result reported in~\cite{Grigoriev97} for the conventional deterministic linear quadratic optimal control.

\subsection{Chaotic Map Lattice}
Generally chaotic systems are more difficult to synchronize than non--chaotic systems. For the coupled map lattice defined in Equation~\eqref{eq:eq4} chaos can be obtained when $3.5 < a \le 4$. In this section we show the control result in this chaotic regime by taking $a=4$ and $\epsilon = 0.25$. The length of the lattice is taken to be $10$, $L=10$, and the covariance matrix of the controller is assumed to be $\Gamma = 0.01 I_{2 \times 2}$.

\begin{figure}[htbp]
  \centerline{
    \begin{tabular}{cc}
     \resizebox{3in}{2in}{\includegraphics{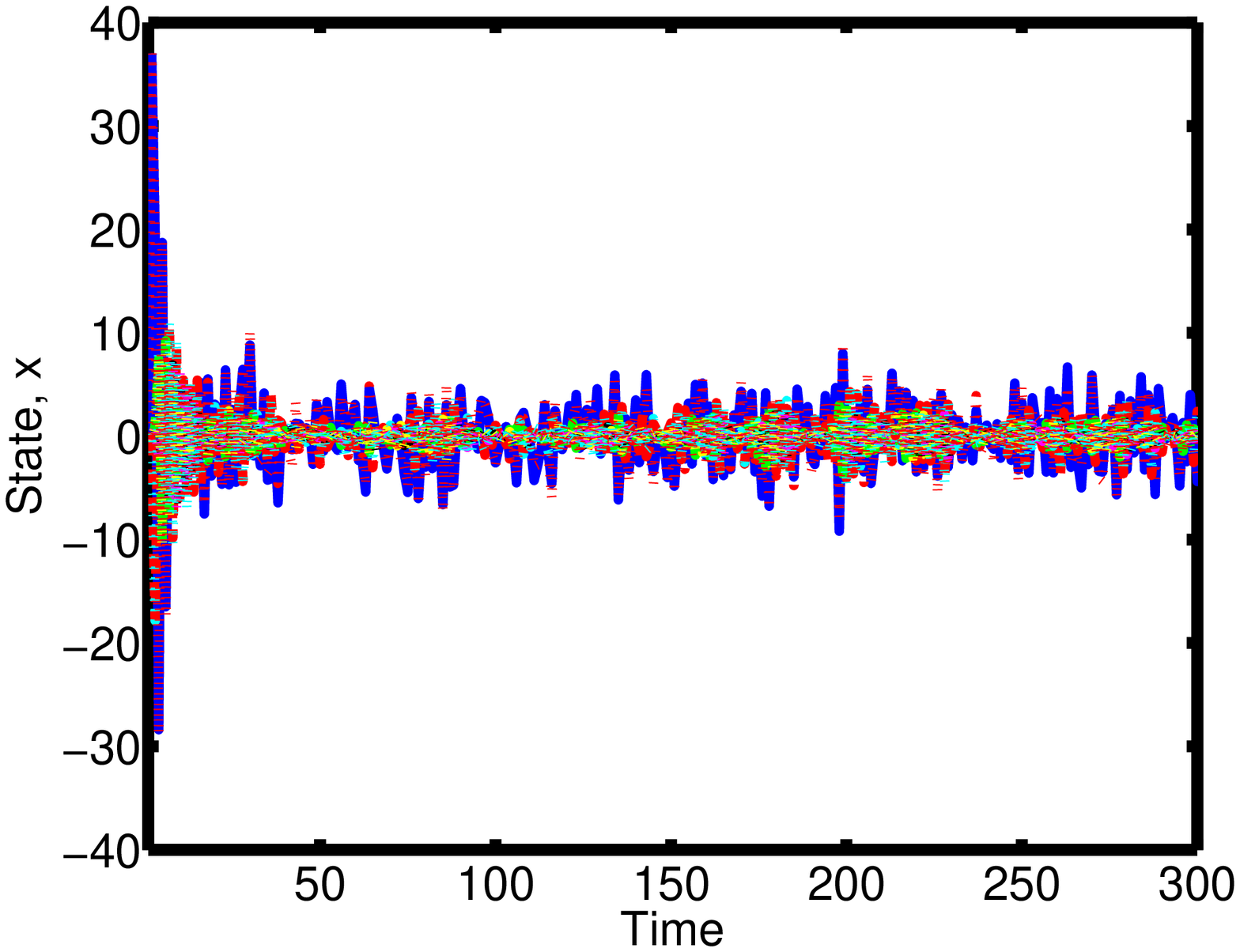}}
     &
      \resizebox{3in}{2in}{\includegraphics{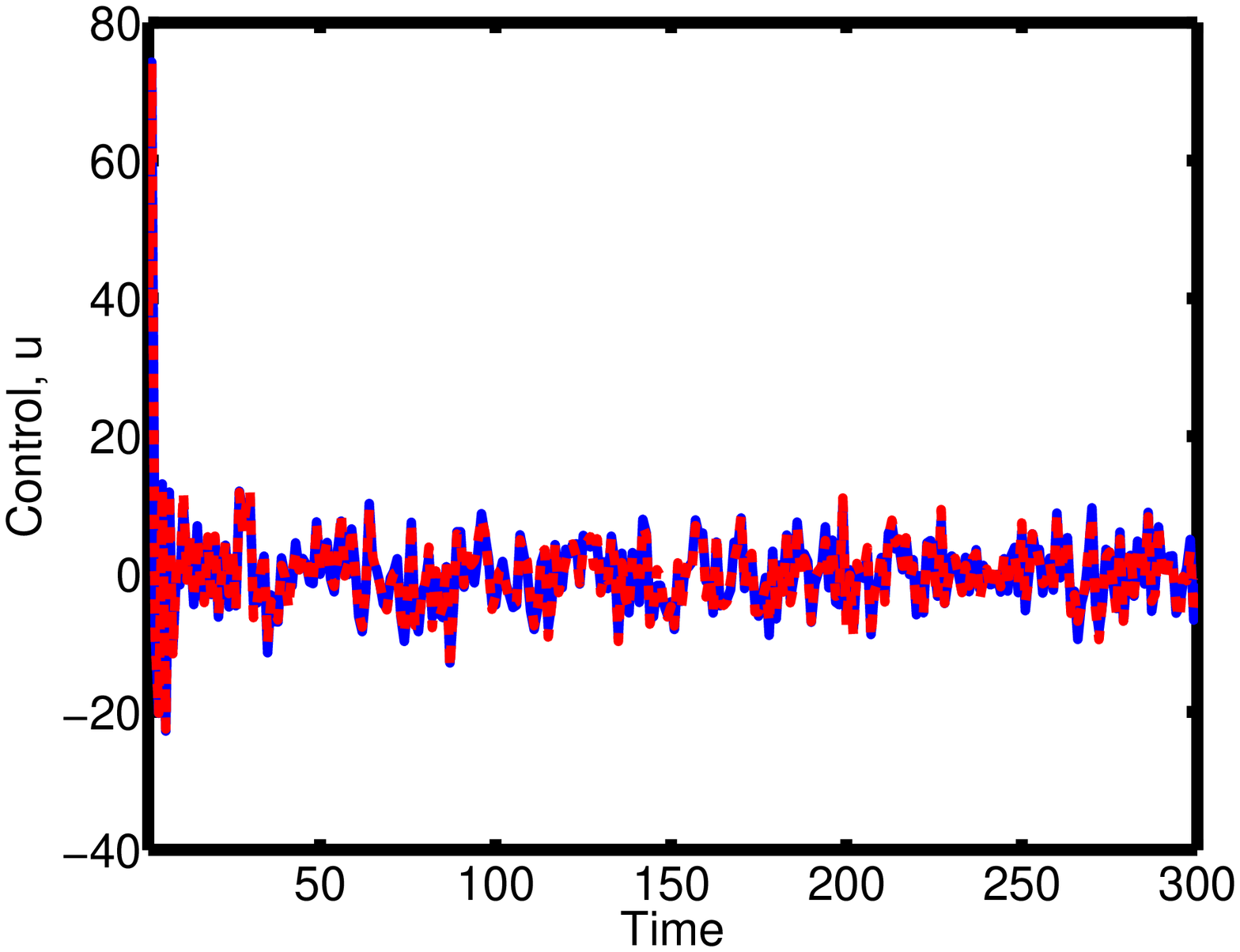}}
      \\
       (a) & (b)
      \\
     \resizebox{3in}{2in}{\includegraphics{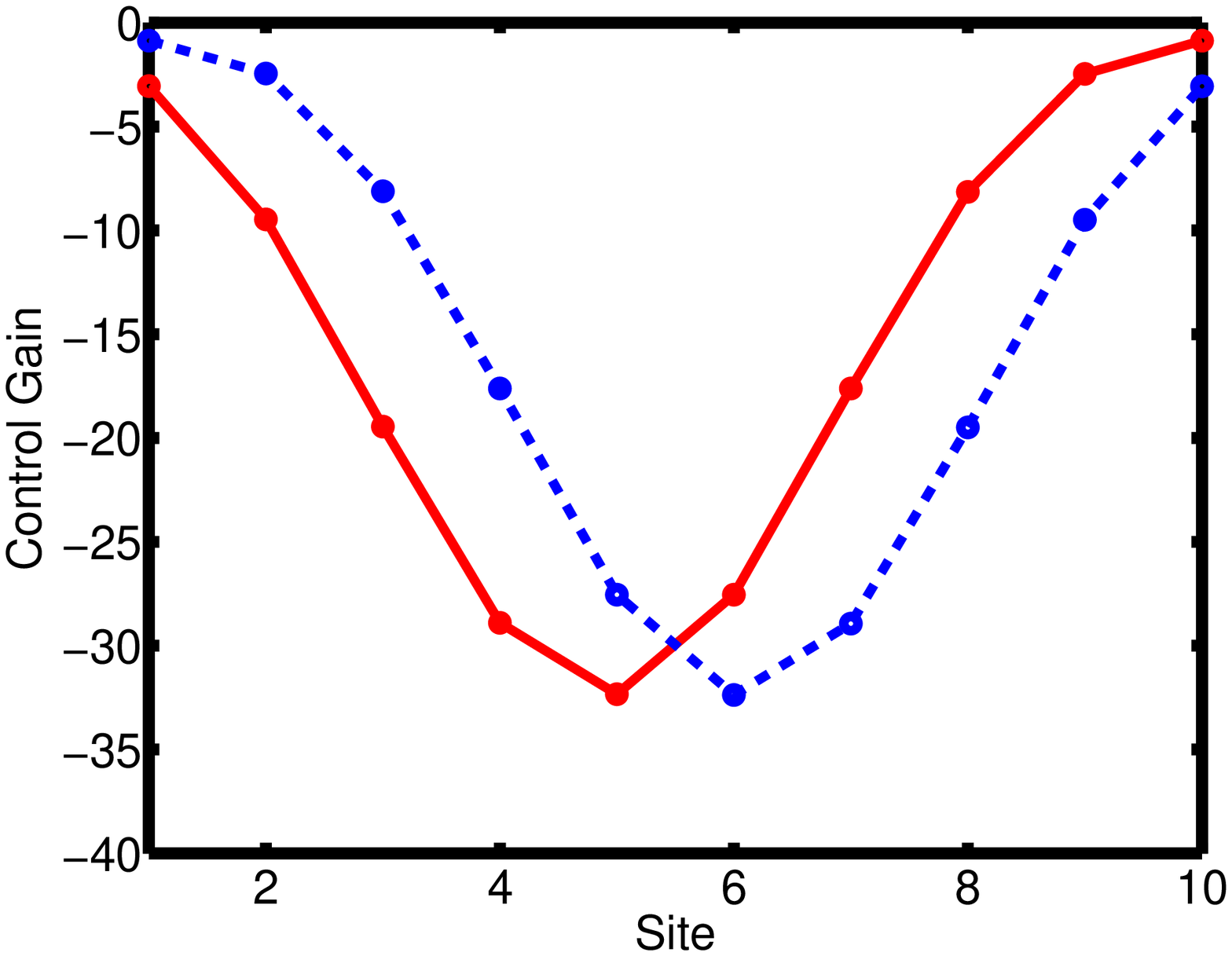}}
     &
     \resizebox{3in}{2in}{\includegraphics{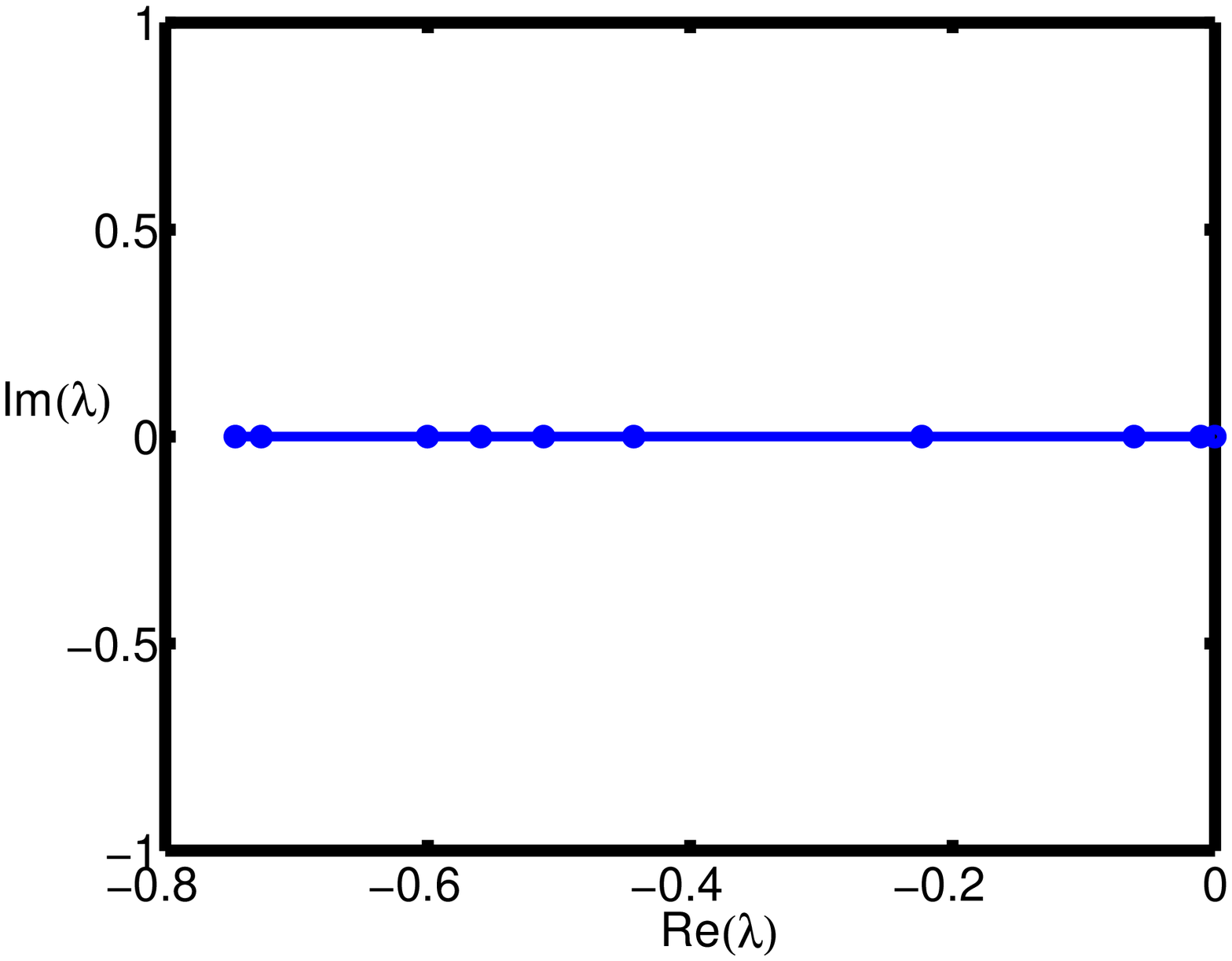}}
     \\
      (c) & (d)
     \end{tabular}
    }
  \caption[Chaotic coupled map lattice.] {Chaotic coupled map lattice with, $L=10$, $a=4$, and $\epsilon=0.25$: (a) Evolution of states. (b) Control efforts. (c) Control Gain. (d) Eigenvalues of the asymptotic stability matrix.}
  \label{fig:Fig3}
\end{figure}

The probabilistic model of the lattice is estimated using GLNN and a global diagonal covariance matrix $\Sigma = [0.00094, 0.0017, 0.00109, 0.00089, 0.0008, 0.00097, 0.00104, 0.00092, 0.0011, \\0.00101]$. Using this matrix the rank of the asymptotic stability matrix, $\text{rank}(AS[((\Sigma^{-1})^{1/2})^T, A])$ is found to be $10$, hence satisfying Assumption$1$. This is further verified by finding the eigenvalues, $\lambda_k$ of the matrix $(A +BC)$ which all are found to be inside the unit circle as can be seen from Figure~\ref{fig:Fig3}(d). Figures~\ref{fig:Fig3}(a) and (b) confirm the result obtained for the non-chaotic map lattice. The designed probabilistic control approach is capable of globally synchronizing the stochastic network. The solution to the control matrix $C$ is shown in Figure~\ref{fig:Fig3}(c). The result is again consistent with what is obtained for the non-chaotic example.

\section{Conclusion}\label{sec:Con}
In this paper pinning synchronization of a stochastic class of complex dynamical coupled map lattice with spatiotemporal chaos network has been investigated in detail. We have proven that a probabilistic adaptive control approach can be successfully used to globally synchronize a network of complex stochastic systems in the presence of uncertainties and noise. We also have defined the concept of stochastic pinning-controllability with the help of the covariance matrix of the residual error of the lattice state. The theoretical findings presented was then validated on two examples in the chaotic and non chaotic regimes of the network. Both examples confirm the effectiveness of the proposed probabilistic control approach in globally synchronizing the states of the network. Current work is addressing the analytical estimation of the bounds of the maximum controllable length and the stochastic controllability of the non--stationary lattice networks.
\\
\\
\textbf{Acknowledgment}: The author would like to thank Prof. D. Lowe for various insightful discussions and pointers to the literature during the author's academic visit to Aston University in the year 2010/2011. The author would also like to thank the editor and anonymous reviewers for comments that helped improve the manuscript.


\end{document}